\documentclass[12pt]{article}
\usepackage{amsmath, amssymb}

\numberwithin{equation}{section}

\newcommand{\R}{\mathbb R}

\newcommand{\N}{\mathbb N}

\newcommand{\QED}{\hspace{.2in}\square\newline}

\newtheorem{definition}{Definition}[section]

\newtheorem{claim}{Claim}[section]

\begin{document}

\begin{center}
{\Large \textbf{A Note on Counting Lattice Points \\ in Bounded Domains}} \vskip 2em

{ J. LaChapelle}\\
\vskip 1em
\end{center}

\begin{abstract}
Zeros and poles of $k$-tuple zeta functions, that are defined here implicitly, enable localization onto prime-power $k$-tuples in pair-wise coprime $k$-lattices $\mathfrak{N}_k$. As such, the set of all $\mathfrak{N}_k$ along with their associated zeta functions encode the positive natural numbers $\N_{>1}$. Consequently, counting points of $\mathbb{Z}_{\geq0}$ can be implemented in $\{\mathfrak{N}_k\}$.  Exploiting this observation, we derive explicit formulae for counting prime-power $k$-tuples and use them to count lattice points in well-behaved bounded regions in $\R^2$. In particular, we count the lattice points contained in the circle $S^1$. The counting readily extends to well-behaved bounded regions in $\R^n$.
\end{abstract}

\section{Introduction}
Constructing explicit summatory functions on $\mathbb{Z}^n$ is mostly trivial with Fourier techniques. But as soon as a boundary is imposed, identifying lattice points that are very close to the boundary becomes difficult. The difficulty stems from Fourier duality, and it is the salient feature in lattice-point counting problems: At issue is the fact that discerning \emph{geometric} localization through increased resolution via Fourier transform produces the opposite effect in the dual space, and the presence of a proximate boundary makes it increasingly difficult to control error terms in the transforms of counting functions.

In this note we propose a method of locating lattice points near boundaries that can be characterized as \emph{algebraic} localization. It is based on two observations: first, the natural numbers $\N_{>1}$ are represented by special points comprised of prime-power $k$-tuples in the set of all pair-wise coprime $k$-lattices $\{\mathfrak{N}_k\}$; and second, one can define $k$-lattice analogs of the Riemann zeta function that implement prime-power $k$-tuple counting. Together, these allow one to count integers in a bounded domain on $\mathbb{Z}_{\geq0}$, and the $k$-tuple zeta zeros (which unfortunately we do not know yet) dictate the achievable error control. In particular, if we are lucky enough that the zeros obey the Riemann hypothesis(RH), then the error can be anticipated to be of the same order as Chebyshev $\psi(x)$ under the RH.

As this is a short note, the presentation will be fairly sparse. We begin by constructing an explicit formula for counting prime powers up to a cut-off $x$. In order to count all integers, it is clear that the prime-power counting function will have to be generalized to prime-power $k$-tuples. Counting \emph{admissible} prime-power $k$-tuples was studied in \cite{LA2}. Both exact and explicit formulae were developed. Here we use these to construct explicit formulae for counting prime-power $k$-tuples --- admissible or otherwise. Unfortunately, without the $k$-tuple zeta zeros, we can't reduce the explicit formulae to sums of residues. Nevertheless, the explicit formulae can be used to formally count lattice points in well-behaved bounded domains in $\R^n$. As an example we count the number of lattice points contained in the circle $S^1$ of radius $R$ and posit the error goes like  $O(R^{1/2+\epsilon})$ if the $k$-tuple zeta zeros abide by the RH.

\section{Counting prime powers}
Start with the well-known explicit formula for counting weighted prime powers;
\begin{equation}\label{num of weighted prime powers}
 J(x)=\sum_{p^k\leq x}\frac{1}{k}
 =\sum_{n=2}^{ x}\frac{\Lambda(n)}{\log(n)}
 =-\lim_{T\rightarrow\infty}\frac{1}{2\pi i}\int_{c-iT}^{c+iT}\mathrm{Ei}\left(\log( x^s)\right)
 \,d\log(\zeta(s))\;\,\;\;\;\;\;c>1
\end{equation}
with $x\in\R_+$ and the principle value prescription for the exponential integral.(see e.g. \cite{MO}) Our goal is to sum over prime powers without the $1/k$ weight. To that end, consider the cousin of the Riemann zeta function
\begin{equation}\label{prime zeta}
\log\left(\mathfrak{z}(s)\right):=-\sum_{n=1}^\infty\frac{\mu(n)\Lambda(n)}{\log(n)n^{s}}\;,\;\;\;\;\Re(s)>1
\end{equation}
where the convergence for $\Re(s)>1$ follows by comparison with $\zeta(s)$. Clearly the right-hand side is just the prime zeta function $P(s)$ (up to a minus sign). Consequently, since $P(s)=\sum_{m=1}^\infty\tfrac{\mu(m)}{m}\log(\zeta(m s))$, equation (\ref{num of weighted prime powers}) yields an explicit formula for primes;\footnote{This is obviously equivalent to the standard Moebius inversion of $J(x)$.}
\begin{eqnarray}\label{3}
 \pi( x)&=&\sum_{p\leq x}1\notag\\
 &=& -\sum_{n=2}^{ x}\frac{\mu(n)\Lambda(n)}{\log(n)}\notag\\
 &=&-\lim_{T\rightarrow\infty}\frac{1}{2\pi i}\int_{c-iT}^{c+iT}\mathrm{Ei}\left(\log( x^s)\right)
 \,d\log(\mathfrak{z}(s))\;\,\;\;\;\;\;c>1\notag\\
 &=&\sum_{m=1}^{\log_2(x)}\frac{\mu(m)}{m}\left\{\mathrm{Ei}(\log(x^{\frac{1}{m}}))
 -\sum_\rho\mathrm{Ei}(\log( x^{\frac{\rho}{m}}))-\log(2)
 -\sum_{k=1}^\infty\mathrm{Ei}(\log( x^{\frac{-2k}{m}}))\right\}\notag\\
\end{eqnarray}
This in turn yields an explicit expression for the total number of prime powers (without the $1/k$ weight) up to a cut-off $x\in\R_+$;
\begin{eqnarray}
\Pi(x)&=&\sum_{p^k\leq x}1\notag\\
&=&\sum_{n=1}^{\log_2(x)}\pi(x^{1/n})\notag\\
&=&\sum_{m,n=1}^{\log_2(x)}\frac{\mu(m)}{m}\left\{\mathrm{Ei}(\log(x^{\frac{1}{mn}}))
-\sum_\rho\mathrm{Ei}(\log( x^{\frac{\rho}{mn}}))-\log(2)
 -\sum_{k=1}^\infty\mathrm{Ei}(\log( x^{\frac{-2k}{mn}}))\right\}\notag\\
 &=&\sum_{m,n=1}^{\infty}\frac{\mu(m)}{m}\left\{\mathrm{Ei}(\log(x^{\frac{1}{mn}}))
-\sum_\rho\mathrm{Ei}(\log( x^{\frac{\rho}{mn}}))-\log(2)
 -\sum_{k=1}^\infty\mathrm{Ei}(\log( x^{\frac{-2k}{mn}}))\right\}\;.\notag\\
\end{eqnarray}

\section{Counting prime-power $k$-tuples}
Now extend the prime-power counting from the previous section to the pair-wise coprime $k$-lattice $\mathfrak{N}_k$. Given a set\footnote{We do not require $\mathcal{H}_k$ to be an admissible set, because we are only counting up to some cut-off $x$.} $\mathcal{H}_k=\{0,h_2,\ldots,h_k\}$ with $h_2<h_3<\cdots<h_k$ and $h_i\in\mathbb{Z}_+$, define the prime $k$-tuple log-zeta function
\begin{definition}
\begin{equation}
\log\left(\mathfrak{z}_{(k)}(s)\right)
:=(-1)^k\sum_{n=1}^\infty\mu(n)\cdots\mu(n+h_k)\frac{\Lambda_{(k)}(n)}{\log_{(k)}(n)\,n_{(k)}^s}\;.
\end{equation}
\end{definition}
where
\begin{eqnarray}
n^s_{(k)}&:=&\left(n(n+h_2)\cdots (n+h_k)\right)^{s/k}\notag\\
\Lambda_{(k)}(n)&:=&\Lambda(n)\cdots\Lambda(n+h_k)\notag\\
\log_{(k)}(n)&:=&\log(n)\cdots\log(n+h_k)\;.
\end{eqnarray}
Notice that $n_{(k)}$ is the geometric mean of $n_k:=(n,n+h_2,\ldots,n+h_k)$ with $n\in\mathbb{N}_+$.

It is useful to also define the prime-power $k$-tuple log-zeta function:
\begin{equation}
\log\left(\zeta_{(k)}(s)\right)
:=\sum_{n=1}^\infty\frac{\Lambda_{(k)}(n)}{\log_{(k)}(n)\,n_{(k)}^s}
\end{equation}
Although $\mathfrak{z}_{(k)}(s)$ and $\zeta_{(k)}(s)$ are related through Moebius, we view $\zeta_{(k)}(s)$ as more fundamental since it characterizes prime powers in $\mathfrak{N}_k$. And prime powers --- not just primes --- seem to follow a gamma distribution\cite{LA2}.

Just like the $k=1$ case from the previous section, the $k$-tuple log-zeta functions can be used to locate prime and prime-power $k$-tuples:
\begin{claim}\label{explicit formula}
Let $\mathrm{P}_k\subset\mathfrak{N}_k$ be the set of prime $k$-tuples where $\mathrm{p}_k=\left(p,p+h_2,\ldots,p+h_k\right)$ with $p$ a prime and $\mathcal{H}_k:=\{0,h_2,\ldots,h_k\}$. Put $\widetilde{x}= x+\epsilon$ with $x\in\N_+$ and $0<\epsilon<1$. Let $\sigma_a$ be the abscissa of absolute convergence of $\sum_{n=1}^\infty\frac{\mu_{(k)}(n)\Lambda_{(k)}(n)}{\log_{(k)}(n)\,n_{(k)}^s}$. Then, for $c>\sigma_a$, the number of prime $k$-tuples up to a cut-off $x$  associated with $\mathcal{H}_k$ is
\begin{eqnarray}\label{explicit eq}
 \pi_{(k)}(x)
 &:=&\sum_{\mathrm{p}_k\in\mathrm{P}_k\,;\,\,p\leq x}1\notag\\
 &=&\lim_{\epsilon\rightarrow0}\lim_{T\rightarrow\infty}\frac{(-1)^{k-1}}{2\pi i}\mathcal{R}_{(k)}(1)\int_{c-iT}^{c+iT}\mathrm{Ei}_{(k)}(\log(\widetilde{x}^{s}))
 \;d\log(\mathfrak{z}_{(k)}(s))\,,
 \;\,\;\;\;\;\;c>\sigma_a\notag\\
 &=&\sum_{n=2}^x\frac{\mu_{(k)}(n)\Lambda_{(k)}(n)}{\log_{(k)}(n)}
 \end{eqnarray}
where
\begin{equation}
\mathrm{Ei}_{(k)}(\log(x)):=\int_2^x\frac{\log_{k-1}(r_{(k)})}{\log_{(k)}(r)}\;dr
\end{equation}
and $\mathcal{R}_{(k)}(x):=\log_{(k)}(x)/\left(\log(x)\log^{k-1}(x_{(k)})\right)$.
\end{claim}
\emph{Proof sketch}: Integrate (\ref{explicit eq}) by parts. The boundary terms won't contribute because:
1) A comparison test between the series representations of $\log^{(k-1)'}(\zeta_{(k)}(s))$ and $\log(\zeta(s))$ yields a finite $\sigma_a$. So, for $s=c+it$ with $c\in\R$, we have $\lim_{t\rightarrow\infty}|\log^{(k-1)'}(\zeta_{(k)}(c+it))|<\infty$ provided $c>k$. To see this, observe that
\begin{eqnarray*}
 &&|\log^{(k-1)'}(\zeta_{(k)}(c+it))|
 \leq\sum_{p^m}\left|\frac{\log^{k-1}(p_{(k)}^m)}{m^k\,p_{(k)}^{ms}}\right|
=\sum_{p^m}\frac{\log^{k-1}(p_{(k)})}{m\,p_{(k)}^{mc}}\notag\\
&&\hspace{1.6in}=\sum_{p^m}\frac{1}{m\,p_{(k)}^{m(c-(k-1))}}
\frac{\log^{k-1}(p_{(k)})}{p^{m(k-1)}_{(k)}}
<\sum_{p^m}\frac{1}{m\,{(p_{(k)}^m)}^{(c-(k-1))}}\notag\\
&&\hspace{1.6in}<\sum_{p^m}\frac{1}{m\,{(p^m)}^{(c-(k-1))}}
=|\log(\zeta((c-(k-1))+it))|\;.
 \end{eqnarray*}
2) The inequality $\log^{k-1}(x_{(k)})/\log_{(k)}(x)\leq1/\log(x)$ implies $\mathrm{Li}_{(k)}(x)\leq \mathrm{li}(x)$. Hence, $\lim_{t\rightarrow\infty}|\mathrm{Li}_{(k)}(x^{s})|=0$ because
\begin{eqnarray*}
\lim_{t\rightarrow\infty}\left|\mathrm{Li}_{(k)}(x^{s})\right|
\leq\lim_{t\rightarrow\infty}\left|\mathrm{li}({x}^{(c+it)})\right|
&=&\lim_{t\rightarrow\infty}\left|\frac{x^{(c+it)}}{(c+it)\log(x)}
\left(1+O\left(\frac{1}{(c+it)\log(x)}\right)\right) \right|\notag\\
&\leq&\frac{x^c}{\log(x)}\lim_{t\rightarrow\infty}
\left|\frac{1}{(c+it)}
\left(1+O\left(\frac{1}{(c+it)}\right)\right) \right|=0\;.\notag\\
\end{eqnarray*}

Now use the truncating integral
\begin{equation*}
\frac{1}{2\pi i}\int_{c-i T}^{c+i T}\left(\frac{x}{n_{(k)}}\right)^s\log(x)\frac{\log^{k-1}(x^s_{(k)})}{\log_{(k)}(x^s)}ds
=\left\{\begin{array}{l}
\mathcal{R}_{(k)}(1)^{-1}
+O\left(\frac{\left(\frac{x}{n}\right)^c}{T\log(x/n)}\right)\;\;\;\;\frac{x}{n}>1 \\
O\left(\frac{\left(\frac{x}{n}\right)^c}{T\log(x/n)}\right)\;\;\;\;0<\frac{x}{n}<1
\end{array}\right.\;.
\end{equation*}
For $c>\sigma_a$ this gives
\begin{eqnarray}
 \lim_{\epsilon\rightarrow0}\lim_{T\rightarrow\infty}\frac{(-1)^{k-1}}{2\pi i}\mathcal{R}_{(k)}(1)\int_{c-iT}^{c+iT}\log(\mathfrak{z}_{(k)}(s))
\widetilde{x}^{s}\log(\widetilde{x})
\frac{\log^{k-1}(\widetilde{x}_{(k)}^s)}{\log_{(k)}(\widetilde{x}^s)}\,ds\notag\\
&&\hspace{-4in}=\lim_{\epsilon\rightarrow0}\lim_{T\rightarrow\infty}\frac{\mathcal{R}_{(k)}(1)}{2\pi i}\int_{c-iT}^{c+iT}\sum_{n=1}^\infty\frac{\mu_{(k)}(n)\Lambda_{(k)}(n)}{\log_{(k)}(n)\,n_{(k)}^s}
 \widetilde{x}^{s}\log(\widetilde{x})
\frac{\log^{k-1}(\widetilde{x}_{(k)}^s)}{\log_{(k)}(\widetilde{x}^s)}\,ds\notag\\
 &&\hspace{-4in}=\lim_{\epsilon\rightarrow0}\lim_{T\rightarrow\infty}
 \sum_{n=1}^\infty\frac{\mu_{(k)}(n)\Lambda_{(k)}(n)}{\log_{(k)}(n)}
 \frac{\mathcal{R}_{(k)}(1)}{2\pi i}\int_{c-iT}^{c+iT}
 \frac{\widetilde{x}^{s}}{n^s_{(k)}}\log(\widetilde{x})
\frac{\log^{k-1}(\widetilde{x}_{(k)}^s)}{\log_{(k)}(\widetilde{x}^s)}\,ds\notag\\
 &&\hspace{-4in}=\lim_{\epsilon\rightarrow0}\sum_{n\leq\lfloor\widetilde{r}\rfloor}
 \frac{\mu_{(k)}(n)\Lambda_{(k)}(n)}{\log_{(k)}(n)}=\sum_{n\leq x}\frac{\mu_{(k)}(n)\Lambda_{(k)}(n)}{\log_{(k)}(n)}
 \end{eqnarray}
where the third equality follows from the lemma. (Justifying the interchange of the sum and integral is straightforward, and interchange of the $T$-limit and sum is allowed because the summand contains $O(n^{-c})$ with $c>\sigma_a$.)
$\QED$

Let us write ${x}^{\frac{1}{\mathbf{m}}}=(x^{\frac{1}{ m_1}}(x+h_2)^{\frac{1}{m_2}}\cdots(x+h_k)^{\frac{1}{m_k}})^{1/k}$. Then determine that
\begin{equation}
\Pi_{(k)}(x)=\sum_{|\mathbf{m}|=k}^\infty\pi_{(k)}({x}^{\frac{1}{\mathbf{m}}})
\end{equation}
counts the number of prime-power $k$-tuples along the ray determined by $\mathcal{H}_k$ up to the cut-off $x$. Note the sum upper limit can be replaced by some $0<M(x,k)<\infty$. By construction, the prime-power $k$-tuples counted by the sum are \emph{ordered} tuples $(p_1^{\alpha_1},p_2^{\alpha_1},\ldots,p_k^{\alpha_k})$ such that $p_1<p_2<,\ldots,p_k$ so the fundamental theorem of arithmetic implies each ordered tuple corresponds to a counting number in $\N_{>1}$. 

Because $\mathcal{H}_k$ determines a ray in $\mathfrak{N}_k$ and $\pi_{(k)}(x)$ only counts the number of prime $k$-tuples along that ray, to get the total number of prime tuples requires a sum over all combinations of  $\mathcal{H}_k$ allowed by the cut-off and hence\footnote{The $+1$ term on the left-hand side accounts for the fact that the counting starts at the natural number $1$ which is not a prime.}
\begin{claim}
\begin{equation}\label{localization}
1+\sum_{ k=1 \,;\,\mathcal{H}_k}^\infty\Pi_{(k)}(x)
=\lfloor x\rfloor
\end{equation}
where the sum over $k$ includes all allowed $\mathcal{H}_k$ with $\left(p^{\alpha_1}(p+h_2)^{\alpha_2}\cdots (p+h_k)^{\alpha_k}\right)\leq x$ for any prime $p$. The upper limit of the sums can be replaced by some $0<M(x)<\infty$.
\end{claim}

It is interesting to note that the average number of prime-power $k$-tuples associated with a particular $\mathcal{H}_k$ is conjectured\cite{LA2} to be
\begin{equation}
\overline{\Pi_{(k)}(x)}= C_{(k)}\int_2^x\frac{1}{\log_{(k)}(r)}\,dr\sim C_{(k)}O(x/\log^k(x))
\end{equation}
where $C_{(k)}$ is the singular series. Asymptotically, this is the Hardy-Littlewood $k$-tuple conjecture, and it is not immediately obvious that the sum over all allowed $\mathcal{H}_k$ will produce $O(x)$. However, Gallagher\cite{GA} showed that, for fixed $k$, $\sum_{h\leq x}C_{(k)}\sim x^k$ as $x\rightarrow\infty$  where $h:=\mathrm{max}(h_i\in\mathcal{H}_k)$. In this limit we have
\begin{equation}
\sum_{h\leq x}\overline{\Pi_{(k)}(x)}\sim \int_2^x\frac{x^k}{\log_{(k)}(r)}\,dr\sim \,O(x^{k+1}/\log^k(x))\;,
\end{equation}
and the sum over $k$ yields $\sum_{k=1}^\infty\sum_{h\leq x}\overline{\Pi_{(k)}(x)}\sim O(x)$ --- so all is well.

At the risk of pointing out the obvious, {if} the $k$-tuple zeta functions exist and are well defined, (\ref{localization}) furnishes an explicit lattice-point counting function. Additionally, {if} the $k$-tuple zeta zeros follow the RH, the fact that the averages combine to give $O(x)$ suggests the sums contributing to $\lfloor x\rfloor$ will yield an error of $O(x^{1/2+\epsilon})$. In particular, the conjectured error term for the number of lattice points within the circle $S^1\subset\R^{2}$ of radius $R$ will go like $O\left(R^{1/2+\epsilon}\right)$. This coincides with previous conjectured and expected sharpest error in the literature.(see e.g. \cite{HB, IKKN})

\section{Counting lattice points}
Let $f:[0,x]\subset\R_+\rightarrow\R_+$ be a piece-wise continuous bijection. We want to count integer lattice points under the graph of $f$ in $\mathbb{Z}_{>0}^2$ (i.e. the $2\mathrm{d}$ integer lattice excluding the $x$ and $y$ axes). Given (\ref{localization}), the task is straightforward:
\begin{claim}
\begin{equation}\label{lattice points}
\#\,\mathrm{of\,points}=\sum_{n=1}^{x}\lfloor f(n)\rfloor
=\sum_{n=1}^{x}\left(1+\sum_{k=1 \,;\,\mathcal{H}_k}^{M(n)}\Pi_{(k)}(f(n))\right)
=\lfloor x\rfloor+\sum_{n=1}^{x}\;\sum_{k=1 \,;\,\mathcal{H}_k}^{M(n)}\Pi_{(k)}(f(n))
\end{equation}
where $0<M(n)<\infty$ is determined by the (finite) cut-off $f(n)$. The extension to higher dimension lattices is immediate.
\end{claim}
Of course the right-hand side includes a rather involved sum over allowed $\mathcal{H}_k$ including evaluations of the associated explicit integrals --- which presents a totally impractical algorithm for calculating the floor function: All we are doing is expressing the lattice points less than or equal to some cut-off as prime-power $k$-tuples and then counting them. Its advantage lies in the analytic expression being ultimately defined in terms of $k$-tuple zeta zeros. This may afford better control over asymptotics than harmonic counting methods.

As a simple illustration, count the lattice points in the upper quadrant of the circle $f(r)=\sqrt{R^2-r^2}$ excluding points on the $x$ and $y$ axes. More appositely, estimate the error $E(R)$ where
\begin{equation}\label{error}
\#\,\mathrm{of\,points}+E(R)=\int_0^R f(r)\,dr-\lfloor R\rfloor\;.
\end{equation}
Hence,
\begin{eqnarray}
E(R)=\tfrac{\pi R^2}{4}-\lfloor R\rfloor-\sum_{n=1}^{R}
\lfloor f(n)\rfloor
=\tfrac{\pi R^2}{4}-2\lfloor R\rfloor-\sum_{n=1}^{R}
\;\sum_{k=1 \,;\,\mathcal{H}_k}^{M(n)}\Pi_{(k)}(\sqrt{R^2-n^2}\,)\;.
\end{eqnarray}
In particular, for $R=6.8$ we have $M(n)=2$ for $n\in\{1,2,3\}$  and the only participating prime-power $2$-tuple is $(2,3)$ while $M(n)=1$ for $n\in\{4,5,6\}$. The prime-power $2^2$ is included for all $n$. The rest of the contributions to the inner sum are single primes. Adding everything up gives $E(6.8)=36.3168-12-24$.

To recap, for algebraic localization the floor function is represented as a finite sum; in contrast to the Fourier representation $\lfloor x\rfloor=x-\frac{1}{2}+\frac{1}{\pi}\sum_{k=1}^\infty\sin(2\pi k x)/k$. Moreover, $\Pi_{(k)}(x)$ has an explicit integral representation in terms of $k$-tuple zeta functions which can be used to bound the error assuming $\mathfrak{z}_{(k)}(s)$ can be explicitly characterized. If we are lucky and the $k$-tuple zeta zeros are governed by the RH, then this exercise suggests the lattice-point counting error for any well-behaved region in $\R^{2}$ of length unit $R^j$ will go like $O\left(R^{j/2+\epsilon}\right)$. In particular, for the hyperbola $xy=R$ this suggests $O\left(R^{1/4+\epsilon}\right)$. Of course our formal and elementary arguments are no substitute for a rigorous and thorough treatment, but hopefully the proposed strategy of algebraic localization will serve as motivation to do so.

\end{document}